\documentclass[11pt]{article}
\usepackage{amsmath,amsthm,amsfonts,amssymb}
\usepackage[all]{xy}
\usepackage[usenames]{color} 
\usepackage{makeidx}

\usepackage{xr}

\newcommand\brefGM[1]{\ref{GM#1}}
\newcommand\brefRK[1]{\ref{RK#1}}

\newcommand\brefRI[1]{\ref{RI#1}}

\newcommand\refGM[1]{I.\brefGM{#1}}
\newcommand\refRK[1]{II.\brefRK{#1}}

\newcommand\refRI[1]{IV.\brefRI{#1}}

\newcounter{enumitemp}

\newtheorem*{theorem*}{Theorem}
\newtheorem*{TheoremA}{Theorem A}
\newtheorem*{TheoremB}{Theorem B}

\newtheorem*{TheoremC}{Theorem C}
\newtheorem*{TheoremCPrime}{Theorem C${}^\prime$}
\newtheorem*{TheoremD}{Theorem D}
\newtheorem*{TheoremE}{Theorem E}
\newtheorem*{TheoremF}{Theorem F}
\newtheorem*{TheoremG}{Theorem G}
\newtheorem*{TheoremH}{Theorem H}
\newtheorem*{TheoremI}{Theorem I}
\newtheorem*{TheoremJ}{Theorem J}

\newtheorem*{proposition*}{Proposition}

\newtheorem*{fact*}{Fact}

\theoremstyle{definition}

\newtheorem*{defn*}{Definition}

\theoremstyle{remark}

\newcounter{remarks}
{\paragraph*{Remarks}\ skip
 \begin{list}{\arabic{remarks}. }{\usecounter{remarks}%
 \setlength{\leftmargin}{0in}%
 \setlength{\rightmargin}{0in}%
 \setlength{\labelsep}{0pt}%
 \setlength{\labelwidth}{0pt}%
 \setlength{\listparindent}{0pt}%
 }
}
{
\end{list}
}

\newcommand\from\colon
\newcommand\inv{{-1}}
\newcommand\subgroup{<}
\newcommand\infinity\infty
\newcommand\na{\text{na}}
\newcommand\supp{\text{supp}}
\newcommand\disjunion\amalg
\newcommand\act\curvearrowright
 
\DeclareMathOperator\PG{PG}
\newcommand\PGF{\PG^{\cal F}}

\DeclareMathOperator\IA{IA}

\newcommand{\Z}{{\mathbb Z}}

\newcommand\GL{\mathsf{GL}}
\newcommand{\Out}{\mathsf{Out}}
\newcommand{\Aut}{\mathsf{Aut}}

\newcommand{\F}{\mathcal F}

\renewcommand\L{\mathcal L}

\def\B{\mathcal B}
\newcommand{\A}{\mathcal A}
\newcommand{\h}{\mathcal H}

\newcommand{\eg}{EG}
\newcommand{\noneg}{NEG}

\renewcommand\neg\noneg

\newcommand{\ct}{CT}
\newcommand{\cts}{CTs}

\newcommand\BookOne{\cite{BFH:TitsOne}}
\newcommand\BookTwo{\cite{BFH:TitsTwo}}

\newcommand\recognition{\cite{FeighnHandel:recognition}}

\newcommand\PartOne{Part~I~\cite{HandelMosher:SubgroupsI}}
\newcommand\PartTwo{Part~II~\cite{HandelMosher:SubgroupsII}}
\newcommand\PartThree{Part~III~\cite{HandelMosher:SubgroupsIII}}
\newcommand\PartFour{Part~IV~\cite{HandelMosher:SubgroupsIV}}

\newcommand\bdy\partial

\newcommand\intersect\cap
\newcommand\union\cup
\newcommand\<\langle
\renewcommand\>\rangle
\newcommand\meet\wedge

\newcommand\cross\times
\newcommand\restrict{\bigm |}

\newcommand\inject\hookrightarrow

\DeclareMathOperator\rank{rank}

 \newcommand\surjection\twoheadrightarrow
 
\DeclareMathOperator\MCG{\mathcal{MCG}}

\DeclareMathOperator\gl{GL}

\renewcommand\int{{\text{int}}}

\title{Subgroup decomposition in $\Out(F_n)$:\\ Research announcement and introduction}

\author{Michael Handel and Lee Mosher}

\makeindex

\begin{document}

\maketitle

\begin{abstract}
This is the introduction to a series of four papers that develop a decomposition theory for subgroups of $\Out(F_n)$ which generalizes the theory for elements of $\Out(F_n)$ found in \BookOne\ and \BookTwo\ and which is analogous to the decomposition theory for subgroups of mapping class groups found in \cite{Ivanov:subgroups}. In this introduction we state the main theorems and we outline the contents of the whole series.
\end{abstract}


%


\newcommand\arXiv{arXiv}

Given a finite type oriented surface $S$, Ivanov's theorem \cite{Ivanov:subgroups} says that for every subgroup $\h \subgroup \MCG(S)$ of the mapping class group, either $\h$ has finite order, or $\h$ leaves invariant the isotopy class of some essential closed curve system, or $\h$ contains a pseudo-Anosov mapping class. Applications include the work of Bestvina and Fujiwara \cite{BestvinaFujiwara:bounded} that calculates bounded 2nd cohomology of arbitrary subgroups $\MCG(S)$ and studies homomorphisms to $\MCG(S)$ from lattices in higher rank Lie groups.

In this series of papers \cite{HandelMosher:SubgroupsI,HandelMosher:SubgroupsII,HandelMosher:SubgroupsIII,HandelMosher:SubgroupsIV} we present a general study of finitely generated subgroups of $\Out(F_n)$. This series supplants our preprint \cite{HandelMosher:SubgroupOutF_n}, still on the \arXiv, the main theorem of which is the following:

\begin{TheoremA}
For each finitely generated\,\footnote{The proof of Theorem A given in \cite{HandelMosher:SubgroupOutF_n} uses that $\h$ is finitely generated, but the statement erred in omitting that hypothesis.} subgroup $\h \subgroup \Out(F_n)$, either $\h$ has a finite index subgroup that fixes the conjugacy class of some free factor of $F_n$, or $\h$ contains a fully irreducible element.
\end{TheoremA}
\noindent
As we shall explain, this theorem resolves the ``absolute case'' of a broader decomposition problem for subgroups of $\Out(F_n)$ on which this series of papers is focussed.

Theorem A was formulated as a strong parallel to Ivanov's Theorem, using ``fully irreducible'' outer automorphisms as an analogue of pseudo-Anosov mapping classes. It has found application by Bridson and Wade \cite{BridsonWade:ActionsOnFreeGroups} in their study of homomorphisms to $\Out(F_n)$ of lattices in higher rank Lie groups: for any irreducible lattice in a semisimple Lie group of rank~$\ge 2$, any homomorphism to $\Out(F_n)$ has finite image. 

However, Theorem A has limitations. Ivanov's Theorem applies inductively: if one had already applied it to $\h \subgroup \MCG(S)$ and obtained an invariant curve system $C \subset S$ then one could restrict $\h$ to 
the components of $S-C$ and apply Ivanov's Theorem again to these restrictions. This doesn't work as well with Theorem~A. Suppose we had already applied Theorem~A to a subgroup $\h \subgroup \Out(F_n)$ to obtain a free factor $A \subgroup F_n$ whose conjugacy class is virtually $\h$-invariant. We could then apply Theorem~A to ``induct down into $A$'' by restricting (a finite index subgroup of) $\h$ to $\Out(A)$, either finding a fully irreducible element or producing a free factor $A' \subgroup A$ whose conjugacy class is virtually $\h$-invariant. But that yields limited information: Theorem~A gives no information about how $\h$ behaves outside of $A$, and once $A'$ is identified it gives no information about how $\h$ behave between $A'$ and $A$.


\section*{The main theorem (slightly simplified).} 

Our purpose in supplanting the preprint  \cite{HandelMosher:SubgroupOutF_n} with this new series is to expand the applicability of Theorem~A by relativizing it. Our main result is Theorem~C, which will first be stated in the slightly simplified form of Theorem~C${}^\prime$.


\smallskip
\textbf{Review: Free factor systems (\BookOne\ Section~2.6 or see Section~\refGM{SectionSSAndFFS})\footnote{``Section I.X.Y.Z'' or ``Theorem I.V.W'' refers to Section X.Y.Z or Theorem V.W of \PartOne.}} Throughout this introduction we review material as needed from the general theory for individual outer automorphisms that is developed in works of Bestvina, Feighn, and Handel \BookOne, \recognition. In Section~\refGM{SectionPrelim} of \PartOne\ we give an extended and thorough but terse review of this theory. We begin with free factor systems.

A \emph{free factor system} of $F_n$ is a finite set of the form $\F = \{[A_1],\ldots,[A_k]\}$ such that there is a free factorization $F_n = A_1 * \cdots * A_k * B$ where $A_1,\ldots,A_k$ are nontrivial ($B$ may be trivial), and $[\cdot]$ denotes conjugacy class. The elements of $\F$ are its \emph{components}. Inclusion of subgroups induces a partial ordering on free factor systems denoted $\F \sqsubset \F'$, which we express by saying that $\F$ is \emph{contained in} $\F'$, or that $\F'$ is an \emph{extension} of $\F$, or that \emph{the pair $\F \sqsubset \F'$ is an extension}. Every free factor system $\F$ is \emph{realized} by some marked graph~$G$ and some core subgraph $H \subset G$, meaning that $H$ has no valence~$1$ vertices and one may list its components as $H = C_1 \union \ldots \union C_k$ so that $\F = \{[\pi_1(C_1)],\ldots, [\pi_1(C_m)]\}$. Furthermore, any nested sequence of free factor systems can be simultaneously realized by nested core subgraphs of a single marked graph. Given a proper extension $\F \sqsubset \F'$, if there exists a realization $H \subset H' \subset G$ such that the subgraph $H' \setminus H$ is an edge of $H'$ then $\F \sqsubset \F'$ is a \emph{one edge extension}, otherwise $\F \sqsubset \F'$ is a \emph{multi-edge extension}. 

The group $\Out(F_n)$ acts naturally on free factor systems preserving the various properties described above, including the relation~$\sqsubset$. For any subgroup $\h \subgroup \Out(F_n)$ and any $\h$-invariant extension of free factor systems $\F \sqsubset \F'$, we say that $\h$ is \emph{reducible} relative to this extension if there exists a $\h$-invariant free factor system strictly between $\F$ and~$\F'$, otherwise $\h$ is \emph{irreducible} relative to this extension. If furthermore each finite index subgroup of $\h$ is irreducible relative to the extension $\F \sqsubset \F'$ then we say that $\h$ is \emph{fully irreducible} relative to this extension. These concepts apply to an individual outer automorphism $\phi \in \Out(F_n)$ via the cyclic subgroup $\<\phi\>$. Using relative train track theory one sees that there is an attracting lamination of $\phi$ associated to any $\phi$-invariant multi-edge extension relative to which $\phi$ is fully irreducible (see under \emph{Attracting laminations} below). Also, for $\phi$ to be fully irreducible in the usual absolute sense translates in this language into the statement that $\<\phi\>$ is fully irreducible with respect to the extension $\F = \emptyset \sqsubset \{[F_n]\} = \F'$. 


Here is our slightly simplified version of Theorem~C:

\begin{TheoremCPrime}
For each finitely generated subgroup $\h \subgroup \Out(F_n)$ and each $\h$-invariant multi-edge extension of free factor systems $\F \sqsubset \F'$, if $\h$ is fully irreducible relative to $\F \sqsubset \F'$ then there exists $\phi \in \h$ which is fully irreducible relative to $\F \sqsubset \F'$.
\end{TheoremCPrime}

The full version of Theorem~C will improve on this by specifying exactly which finite index subgroup of $\h$ it is sufficient to consider, namely, its intersection with $\IA_n(\Z/3)$. This will allow us to remove the adjective ``fully'' from the statement, without losing any information.

\section*{Rotationless versus $\IA_n(\Z/3)$ (\PartTwo)} 

Outer automorphisms of $F_n$ can exhibit nontrivial finite permutation behavior in several ways. Given $\phi \in \Out(F_n)$, for example, $\phi$ may permute the components of a $\phi$-invariant free factor system, or there may be a conjugacy class in $F_n$ which is $\phi$-periodic but not fixed by $\phi$. For an individual outer automorphism one may control individual instances of this behavior by passing to a finite power. 


There are two important subsets of $\Out(F_n)$ whose elements already have good control over finite permutations without passing to a further power: rotationless outer automorphisms (\recognition\ Section~3, or see Section~\refGM{SectionPrincipalRotationless}); and elements of the finite index characteristic subgroup $\IA_n(\Z/3) \subgroup \Out(F_n)$ which, by definition, is the kernel of the natural epimorphism $\Out(F_n) \to H_1(F_n;\Z/3) \approx \gl(n,\Z/3)$. For instance, members of those two subsets satisfy the following invariance properties:
\begin{itemize}
\item The action of $\phi$ fixes each element of the finite set $\L(\phi)$ of attracting laminations.
\item Every $\phi$-periodic conjugacy class in $F_n$ is fixed by~$\phi$.
\item Every $\phi$-periodic free factor system is fixed by $\phi$ and its components are fixed by~$\phi$.
\item For any $\phi$-invariant extension by free factor systems, $\phi$ is irreducible with respect to that extension if and only if it is fully irreducible.
\end{itemize}
In particular, for $\phi \in \Out(F_n)$ which is either rotationless or in $\IA_n(\Z/3)$ we have the following fact:
\begin{itemize}
\item \emph{$\phi$ is fully irreducible if and only if it is irreducible}
\end{itemize}
Because of this last fact, we often drop the word ``fully''. 

For rotationless outer automorphisms, the proofs of the above properties are found in \recognition; see Fact~\refGM{FactPeriodicIsFixed} for citations. Rotationless outer automorphisms have the advantage of being represented by the best relative train track maps, the so-called ``completely split relative train track maps'' or \cts\  (see \recognition\ Theorem~4.28 or Theorem~\refGM{TheoremCTExistence}). One major disadvantage is that they are not closed under the group operation on $\Out(F_n)$.

For outer automorphisms in the subgroup $\IA_n(\Z/3)$ the proofs of these invariance properties are a significant portion of \PartTwo:

\begin{TheoremB}[Lemma~\refRK{LemmaFFSComponent} and Theorems~\refRK{ThmPeriodicFreeFactor} and \refRK{ThmPeriodicConjClass}] \quad \\ Each of the above invariance properties holds for each element of the subgroup $\IA_n(\Z/3)$. 
\end{TheoremB}

The subgroup $\IA_n(\Z/3) \subgroup \Out(F_n)$ plays a similar role for us as the role played by the kernel of the homomorphism $\MCG(S) \mapsto \Aut(H_1(S;\Z/3))$ in \cite{Ivanov:subgroups}, namely that of providing a finite index subgroup of elements with good invariance properties.
%
%
Every subgroup of $\Out(F_n)$ has a finite index subgroup in $\IA_n(\Z/3)$, namely its intersection; for this reason we often state our results under the restriction that the subgroup is already in $\IA_n(\Z/3)$. As an application of the Theorem~B, given an element or subgroup of $\IA_n(\Z/3)$, irreducibility relative to an extension of free factor systems $\F \sqsubset \F'$ is equivalent to full irreducibility relative to $\F \sqsubset \F'$.

\section*{The main theorem (full version)}

\begin{TheoremC}
For each finitely generated subgroup $\h \subgroup \IA_n(\Z/3)$ and each $\h$-invariant multi-edge extension of free factor systems $\F \sqsubset \F'$, if $\h$ is irreducible relative to this extension then then there exists $\phi \in \h$ which is irreducible relative to this extension. 
\end{TheoremC}

Here is an equivalent version of Theorem~C. A \emph{filtration by free factor systems} is a strictly increasing sequence $\emptyset =\F_0 \sqsubset \F_1 \sqsubset \cdots \sqsubset \F_\ell = \{[F_n]\}$. Each such filtration can be realized simultaneously by some marked graph~$G$ and some nested sequence of core subgraphs $\emptyset = G_0 \subset G_1 \subset\cdots\subset G_\ell = G$; from this one obtains the length bound $\ell \le 2n-1$. A simple induction shows for any subgroup $\h \subgroup \Out(F_n)$, any $\h$-invariant filtration by free factor systems extends to a maximal such filtration $\emptyset = \F_0 \sqsubset \F_1 \sqsubset \cdots \sqsubset \F_k = \{[F_n]\}$. By maximality, $\h$ is irreducible with respect to each extension $\F_{i-1} \sqsubset \F_i$ of this filtration. Theorem~C is therefore equivalent to the following:

\begin{TheoremD}
For each finitely generated subgroup $\h \subgroup \IA_n(\Z/3)$, each maximal $\h$-invariant filtration by free factor systems $\emptyset = \F_0 \sqsubset \F_1 \sqsubset \cdots \sqsubset \F_m = \{[F_n]\}$, and each $i=1,\ldots,m$ such that $\F_{i-1} \sqsubset \F_i$ is a multi-edge extension , there exists $\phi \in \h$ which is irreducible with respect to $\F_{i-1} \sqsubset \F_i$. 
\end{TheoremD}

\emph{Remark.} We do not know whether, in Theorem~D, there exists a single $\phi \in \h$ which is fully irreducible with respect to each multi-edge extension $\F_{i-1} \sqsubset \F_i$, i.e.\ we do not know whether one can switch the order of quantification of $\F_{i-1} \sqsubset \F_i$ and $\phi$. We conjecture that this can be done. The analogous statement for mapping class groups is indeed true.

\smallskip

\emph{Remark.} We do not know whether the conclusion of Theorem~C holds for all one edge extensions, although a lot can be said and this may be a relatively minor point. One edge extensions $\F \sqsubset \F'$ are of three types, which can be expressed in terms of ranks of components and in terms of realizations in a marked graph. In a \emph{circle} extension, $\F' = \F \union \{[B]\}$ for some rank~1 free factor~$B$: one adds a disjoint circle to a subgraph. In a \emph{barbell} extension, $\F' = (\F - \{[A_1],[A_2]\}) \union \{[B]\}$ for some free factor $B \subgroup F_n$ and nontrivial free factorization $B = A_1 * A_2$ such that $[A_1],[A_2] \in \F$: one attaches the endpoints of an edge to distinct components of a subgraph. In a \emph{handle} extension, $\F' = (\F - \{[A]\}) \union \{[B]\}$ for some free factors $A \subgroup B$ such that $\rank(B)=\rank(A)+1$ and $[A] \in \F$: one attaches the endpoints of an edge to the same component of a subgraph. If $\F \sqsubset \F'$ is a circle or barbell extension then there do not exist any free factor systems strictly between $\F$ and $\F'$, and so any element or subgroup preserving $\F \sqsubset \F'$ is fully irreducible relative to to $\F \sqsubset \F'$. However, if $\F \sqsubset \F'$ is a handle extension then there are infinitely many different free factor systems $\F''$ such that $\F \sqsubset \F'' \sqsubset \F'$; in all cases $\F \sqsubset \F''$ is a circle extension and $\F'' \sqsubset \F'$ is a barbell extension. We do not know the answer to the following question: 
\begin{itemize}
\item If a subgroup $\h \subgroup \IA_n(\Z/3)$ is irreducible relative to a handle extension $\F \sqsubset \F'$ for which $\rank(A) \ge 2$, does there exist an element $\phi \in \h$ which is irreducible relative to $\F \sqsubset \F'$?
\end{itemize}

\medskip

The rest of this introduction will outline the proof of Theorem~C and explain the contents of the other parts of this series, \PartOne, \PartTwo, \PartThree, \PartFour. 

\section*{The relative Kolchin theorem for $\Out(F_n)$ (\PartTwo).} 

The proof of Theorem~C breaks into two major steps: the relative Kolchin theorem of \PartTwo; and the ping-pong arguments of~\PartFour. In \cite{HandelMosher:SubgroupOutF_n}, the proof of Theorem~A breaks into two similar steps, the first of which is simply to cite the (absolute) Kolchin theorem from \BookTwo:

\begin{theorem*}[The Kolchin type theorem for $\Out(F_n)$ \BookTwo] For any finitely generated subgroup $\h \subgroup \Out(F_n)$, if each element of $\h$ is a polynomially growing outer automorphism, and if the image of each element of $\h$ under the natural homomorphism $\Out(F_n) \mapsto \Aut(H_1(F_n;\Z)) \approx \GL(n,\Z)$ is unipotent, then there exists an $\h$-invariant filtration by free factor systems $\emptyset = \F_0 \sqsubset \F_1 \sqsubset \cdots \sqsubset \F_k = \{[F_n]\}$ such that each extension $\F_{i-1} \sqsubset \F_n$ in this filtration is a one edge extension.
\end{theorem*}

Recall that for $\phi$ to be polynomially growing means that for each automorphism $\Phi \in \Aut(F_n)$ representing $\phi$ and each $g \in F_n$, the cyclically reduced word length of $\Phi^i(g)$ is bounded above by a polynomial function of~$i$. In topological terms this is equivalent to saying that for each marked graph $G$ and each conjugacy class $c$ of $F_n$, the length of the circuit in $G$ representing $\phi^i(c)$ is bounded above by a polynomial. The hypothesis of the above theorem is captured by the terminology that each element $\phi \in \h$ is a ``UPG'' element of $\Out(F_n)$, the ``U'' referring to ``unipotent'', and the ``PG'' to polynomially growing. This concept is connected with $\IA_n(\Z/3)$, via the result of \BookTwo\ that each PG outer automorphism of $F_n$ that is contained in $\IA_n(\Z/3)$ is UPG (see \BookTwo\ Proposition~3.5). The above theorem therefore applies to all finitely generated subgroups $\h \subgroup \IA_n(\Z/3)$ consisting solely of polynomially growing outer automorphisms.

In our current, relative setting we need to generalize to subgroups $\h \subgroup \IA_n(\Z/3)$ such that each element of $\h$ leaves invariant a certain free factor system~$\F$, allowing for the possibility that the behavior of $\h$ down in $\F$ may be very complicated, in fact elements of $\h$ may even have exponential growth down in $\F$, but up outside of~$\F$ the growth is controlled in that each element of $\h$ is ``polynomially growing relative to $\F$''. This defines a subset of $\Out(F_n)$ denoted $\PGF$. To formulate this rigorously, we review some concepts of attracting laminations.

\smallskip\textbf{Review: Attracting laminations  (\BookOne\ Section~3, or see Section~\refGM{SectionPrelim}).} Associated to each $\phi \in \Out(F_n)$ there is a finite, $\phi$-invariant set of \emph{attracting laminations} denoted $\L(\phi)$. Attracting laminations have an invariant definition (\BookOne\ Definition~3.1.5, or see Section~\refGM{SectionAttractingLams}) which is expressed---without reference to a choice of relative train track representative---solely in terms of the action of $\phi$ on the ``space of lines'' $\B=\B(F_n)$. Attracting laminations also have an equivalent characterization which is expressed in terms of any relative train track representative of~$\phi$ (\BookOne\ Definition~3.1.12, or see Section~\refGM{SectionLams}), which yields a bijection between $\L(\phi)$ and the \eg\ strata of the relative train track representative (perhaps after passing to a power of $\phi$ that fixes each element of $\L(\phi)$, such as a rotationless power). 

Attracting laminations can be thought of as an asymptotic record of the exponentially growing features of $\phi$. In particular, $\phi \in \Out(F_n)$ is polynomially growing if and only if $\L(\phi) = \emptyset$. 

Associated to each attracting lamination $\Lambda \in \L(\phi)$ is its \emph{free factor support}, the smallest free factor system that ``carries'' each leaf of~$\Lambda$, denoted $\F_\supp(\Lambda)$ (\BookOne\ Section~3.2, or see Section~\refGM{SectionSubgroupsCarryingThings}). Distinct elements of $\L(\phi)$ have distinct free factor supports. For any $\phi \in \Out(F_n)$ there is a natural bijection $\L(\phi) \leftrightarrow \L(\phi^\inv)$ where $\Lambda^+ \leftrightarrow \Lambda^-$ if and only if $\F_\supp(\Lambda^+)=\F_\supp(\Lambda^-)$, in which case we say that the ordered pair $\Lambda^\pm = (\Lambda^+,\Lambda^-)$ is a \emph{dual lamination pair} for $\phi$. The set of dual lamination pairs for $\phi$ is denoted~$\L^\pm(\phi)$.

If $\phi \in \Out(F_n)$ leaves invariant some multi-edge extension of free factor systems $\F \sqsubset \F'$ then, using standard results of relative train track theory, one easily shows the following statement: if $\phi$ is fully irreducible relative to that extension, then there exists an attracting lamination $\Lambda \in \L(\phi)$ which is supported by the free factor system $\F'$ but not by~$\F$. The failure of the converse of this statement is an important point in our discussion below of how Theorem~E is used in the proof of Theorem~C.

\smallskip

Given a free factor system $\F$ of $F_n$, we let $\PGF$ denote the subset of all $\phi \in \Out(F_n)$ such that $\phi(\F)=\F$ and $\phi$ is of polynomial growth relative to $\F$ (see Section~\refRK{SectionPGF}). This concept may be defined in terms of laminations by saying that each $\Lambda \in \L(\phi)$ is supported by the free factor system $\F$, and it has an equivalent and more geometric definition as follows: given a marked graph $G$ with subgraph $G'$ realizing the free factor system~$\F$ we have $\phi \in \PGF$ if and only if for each conjugacy class $[c]$ of $F_n$, letting $\sigma_i$ be the circuit in $G$ realizing $\phi^i[g]$, the number of times that $\sigma_i$ crosses edges of $G \setminus G'$ is bounded above by a polynomial function of~$i$.

The following is the main theorem of \PartTwo, what we call the ``relative Kolchin theorem'':

\begin{TheoremE}[Theorem~\refRK{relKolchin}]
For any finitely generated subgroup $\h \subgroup \IA_n(\Z/3)$ and any free factor system $\F$ such that $\h \subset \PGF$, there exists an $\h$-invariant filtration by free factor systems $\F = \F_0 \sqsubset \F_1 \sqsubset \cdots \sqsubset \F_k=\{[F_n]\}$ such that for each $i=1,\ldots,k$ the extension $\F_{i-1} \sqsubset \F_i$ is a one edge extension.
\end{TheoremE}

The proof of this theorem follows the same outline as the proof of the absolute Kolchin theorem in \BookTwo. However, many new arguments are needed, to replace arguments that work in the absolute case but not in the relative case where exponential growth is allowed as long as it is isolated in~$\F$. 

\medskip

The manner in which we shall apply Theorem~E is to get the proof of Theorem~C off the ground. Suppose we are given a subgroup $\h \subgroup \IA_n(\Z/3)$ which is irreducible relative to an $\h$-invariant multi-edge extension $\F \sqsubset \F'$. Let us focus on the special case that $\F' = \{[F_n]\}$ (see Theorem~I below which explains how to reduce the general case to this special case). Since there exists no $\h$-invariant free factor system strictly between $\F$ and $\{[F_n]\}$, by applying Theorem~E we conclude that $\h \not\subset \PG^{\F}$, and so there exists $\psi \in \h$ that is not of polynomial growth relative to $\F$. We therefore have some $\Lambda \in \L(\psi)$ that is not supported by the free factor system~$\F$. This represents a first step towards the conclusion of Theorem~C, although there is still a lot to do, because the existence of $\Lambda$, while necessary for $\phi$ to be fully irreducible relative to $\F \sqsubset \{[F_n]\}$, is not sufficient.

To make further progress towards the conclusion of Theorem~C, using $\psi$ and $\Lambda$ we shall play a game of ping-pong using the action of $\Out(F_n)$ on the space of lines~$\B$. The basic concepts of this ping-pong game were introduced in \BookOne\ for purposes of proving the Tits alternative. Building the fancier ping-pong table needed for our present purposes will take quite a bit of work. 

\section*{Geometric models (\PartOne).} 

A fully irreducible outer automorphism $\phi \in \Out(F_n)$ is said to be \emph{geometric} if it is represented by a homeomorphism $h \from S \to S$ of a compact surface $S$ with nonempty boundary, relative to some isomorphism $\pi_1(S) \approx F_n$. This definition, introduced in \cite{BestvinaHandel:tt}, is shown there to imply that the surface $S$ has connected boundary and the mapping class of $h$ is pseudo-Anosov. 

In \BookOne\ Definition~5.1.4, the concept of geometricity is extended from a property of a fully irreducible outer automorphism to a property of an \eg\ (exponentially growing) stratum $H_r$ of a relative train track map $f \from G \to G$. The definition of geometricity is expressed in terms of the existence of what we call a ``geometric model''. In the very simplest case where $H_r$ is the top stratum of $G$, this model is obtained by replacing $H_r$ with a surface $S$, attaching all but one component of $\bdy S$ to $G \setminus H_r$, and defining a homotopy equivalence of the resulting 2-complex $X$ that ``extends'' the restriction $f \restrict (G \setminus H_r)$ over the surface $S$ using a pseudo-Anosov homeomorphism of~$S$. 

In \PartOne\ we express geometricity of an \eg\ stratum $H_r$ in terms of a somewhat different geometric model defined by modifying Definition~5.1.4: we carefully embed $H_r$ into $S$ so that the entire 2-complex $X$ deformation retracts onto~$G$; and we modify both the way in which the components of $\bdy S$ are attached to $G \setminus H_r$ and the sense in which the homotopy equivalence of the 2-complex extends $f \restrict (G \setminus H_r)$. As is done in \BookOne, in \PartOne\ we characterize geometricity of an \eg\ stratum in various equivalent ways, one of which is expressed solely in the language of relative train track representatives. 

In any relative train track map $f \from G \to G$ representing a rotationless $\phi \in \Out(F_n)$ there is a bijection between the \eg\ strata $H_r$ of $G$ and the set of attracting laminations of $\phi$. By definition, $H_r$ is geometric if it has a geometric model. In Proposition~\refGM{PropGeomEquiv} we prove that for $\Lambda \in \L(\phi)$ geometricity is a well-defined property of the pair $(\Lambda,\phi)$, independent of the choice of $f \from G \to G$ and the corresponding \eg\ stratum $H_r$. More precisely, if $f' \from G' \to G'$ is another relative train track representative of $\phi$ with \eg\ stratum $H_{r'}$ corresponding to $\Lambda$, then $H_r$ is geometric if and only if $H_{r'}$ is geometric. Even more, geometricity is an invariant of the dual lamination pair $\Lambda^\pm \in \L^\pm(\phi)$, well-defined independent of the choice of a relative train track map representating any positive or negative power of $\phi$ and of the \eg\ stratum corresponding to $\Lambda^\pm$. 

As in \BookOne\ and \BookTwo, the distinction between geometric and nongeometric \eg\ strata is important throughout this series of papers. In addition, we use the geometric models themselves in new ways, by applying transversality arguments in the 2-complex~$X$. In Section~\refGM{SectionGeomModelComplement} we study self-homotopy equivalences of~$X$ with particular focus on \emph{complementary graph of $X$} which is the 1-dimensional subcomplex of $X$ identified with the closure of the complement of the interior of the surface~$S$. In Section~\refGM{SectionPreservingSurface} we prove Lemma~\refGM{LemmaGeomModelHE} which roughly says that any homotopy equivalence of $X$ that preserves the complementary subgraph must also preserve the surface $S$ up to homotopy, inducing a homeomorphism of~$S$; this is related to Waldhausen's theorem that a homotopy equivalence of $S$ that preserves $\bdy S$ must be homotopic to a homeomorphism of $S$. Lemma~\refGM{LemmaGeomModelHE} plays an important role in the proofs in \PartTwo\ of the invariance properties for elements of $\IA_n(\Z/3)$ that are mentioned above, and hence also in the proof of Theorem~E.

\section*{Vertex group systems (\PartOne)}

While free factor systems in $F_n$ are of central importance in the study of $\Out(F_n)$, other kinds of ``subgroup systems'' are also important. A ``subgroup system'' in $F_n$ is just a finite set whose elements are conjugacy classes of finite rank subgroups of $F_n$. For example, it is proved in \cite{GJLL:index} that for any minimal action of $F_n$ on a nondegenerate tree $T$ with trivial arc stabilizers, the stabilizer of every point has finite rank, and there are only finitely many $F_n$-orbits of points whose stabilizers have positive rank; the conjugacy classes of such point stabilizers forms a subgroup system that we call the \emph{vertex group system} of~$T$. Proposition~\refGM{PropVDCC}, with proof suggested to us by Mark Feighn, bounds the length of any nested chain of vertex group systems.

In \PartOne\ we study a certain class of subgroup systems that arise naturally in connection with a geometric model $X$, namely the subgroup systems associated to the complementary subgraph of $X$ and those of its subgraphs that contain $\bdy S$. As a simple example, if $S$ is a compact surface with nonempty boundary and with fundamental group isomorphic to $F_n$ then the components of $\bdy S$ define a vertex group system that is not a free factor system, each of whose components has rank~$1$. We prove in Proposition~\refGM{PropGeomVertGrSys} that the subgroup systems associated in this manner to geometric models are vertex group systems but not free factor systems; this is a key component of the proof of Theorem~F below, and it is important in the ping-pong arguments of \PartFour.

\section*{Weak attraction theory (\PartThree).} Ping-pong arguments in topology are often based on a dynamical theory of attraction and repulsion, applied to the action of some group on some topological space. In topological proofs of the Tits alternative, for example, ping-pong is applied to carefully chosen elements of the group in order to prove that those elements generate a free subgroup. In the context of $\Out(F_n)$ this is carried out in \BookOne. The  Weak Attraction Theorem, Theorem~6.0.1 of~\BookOne, describes attraction--repulsion dynamics for the action of $\Out(F_n)$ on the space of birecurrent lines of~$F_n$, using ``top'' lamination pairs $\Lambda^\pm \in \L^\pm(\phi)$ of the generator $\phi$ as the attractor and repeller.

In \PartThree\ we develop a more expansive weak attraction theory for application to the ping-pong arguments to be carried out in \PartFour. This theory describes attraction--repulsion systems not just on birecurrent lines but on the full space of lines $\B=\B(F_n)$, which is introduced in \BookOne\ and reviewed in Section~\refGM{SectionLineDefs}. Furthermore, it allows any lamination pair $\Lambda^\pm \in \L(\phi)$ for any rotationless $\phi \in \Out(F_n)$ as the attractor--repeller pair, not just a ``topmost'' pair. Given $\phi$, $\Lambda^\pm$, and a line $\ell \in \B$, we characterize when $\ell$ is attracted to $\Lambda^+$ under positive iteration of $\phi$, when to $\Lambda^-$ under negative iteration, and when to neither. This characterization is broken into two theorems: Theorem~F covers just the case of periodic lines; Theorem~G, a more technical statement, considers the general case. 

For application to \PartFour, what is most important is the description and properties of the ``nonattracting subgroup system'' described in Theorem~F, and a less technical statement regarding attraction--repulsion of general lines given in Theorem~H. 

\begin{TheoremF}
For each rotationless $\phi \in \Out(F_n)$ and each dual lamination pair $\Lambda^+ \in \L(\phi)$, $\Lambda^- \in \L(\phi^\inv)$ there exists a subgroup system $\A_\na(\Lambda^\pm)$, called the \emph{nonattracting subgroup system}, with the following properties:
\begin{itemize}
\item\label{TheoremFVGS}
$\A_\na(\Lambda^\pm)$ is a vertex group system.
\item\label{TheoremFCC}
For each conjugacy class $c$ in $F_n$ the following are equivalent:\\
--- $c$ is not weakly attracted to $\Lambda^+$ under iteration of $\phi$;\\
--- $c$ is carried by $\A_\na(\Lambda^\pm)$.
\item $\A_\na(\Lambda^\pm)$ is completely determined by the previous properties.
\item For each conjugacy class $c$ in $F_n$, $c$ is not weakly attracted to $\Lambda^+$ under iteration of $\phi$ if and only if $c$ is not weakly attracted to $\Lambda^-$ under iteration of $\phi^\inv$.
\item\label{TheoremFLams}
The lamination pair $\Lambda^\pm$ is geometric if and only if the vertex group system $\A_\na(\Lambda^\pm)$ is \emph{not} a free factor system. 
\end{itemize}
\end{TheoremF}

The complete description of the vertex group system is $\A_\na(\Lambda^\pm)$ is given in \PartThree, expressed very explicitly in terms of any \ct\ $f \from G \to G$ representing $\phi$; this description generalizes the concept of $\<Z,\rho\>$ developed for the \emph{Weak Attraction Theorem} of \BookOne, Theorem 6.0.1. The full statement of Theorem~F, given in the introduction to \PartThree, contains additional information relating geometricity of the lamination pair $\Lambda^\pm$ to the behavior of the nonattracting subgroup system $\A_\na(\Lambda^\pm)$. The proof of Theorem~F given in \PartThree\ combines relative train track theory with facts about geometric models proved in \PartOne.

Theorem~G, our second weak attraction result, is a complete characterization of those lines $\ell \in \B$ that are \emph{not} weakly attracted to $\Lambda^+$ under iteration of~$\phi$; here we give only a vague statement of Theorem~G. We assume that both $\phi$ and $\phi^\inv$ are rotationless, which implies that for the most important automorphisms representing either of $\phi$ or $\phi^\inv$---namely, the principal automorphisms---every periodic point of $\bdy F_n$ is fixed. There are three obvious sets of lines that are not weakly attracted to~$\Lambda^+$: those carried by the nonattracting subgroup system $\A_\na(\Lambda^\pm)$; the generic leaves of attracting laminations of $\phi^\inv$; and the ``singular lines'' for $\phi^\inv$, which by definition are those lines whose endpoints are nonrepelling fixed points of some principle automorphism representing $\phi^\inv$. The statement of Theorem~G refers to a binary operation of asymptotic concatenation on lines: roughly speaking, two lines which have exactly one end in common can be concatenated at that end to define another line.

\begin{TheoremG}
If $\phi,\phi^\inv$ are both rotationless, then a line $\ell \in \B$ is not weakly attracted to $\Lambda^+$ under iteration of $\phi$ if and only if $\ell$ is the asymptotic concatenation of some finite sequence of lines, each of which is in one of the three ``obvious sets'' described above.
\end{TheoremG}
\noindent
The possible concatenations which can occur in the context of this theorem are quite limited, which leads to a much more precise version of Theorem~G as given in the introduction to \PartThree.

\bigskip

The following result is a distillation of weak attraction theory that avoids most of the technicalities of the statement of Theorem~G while retaining key features needed for application in \PartFour. The proof does depend on all of those technicalities. Item~(2) in the conclusion can be viewed as a uniform version of item~(1).

\begin{TheoremH}
Given rotationless $\phi,\phi^\inv \in \Out(F_n)$ and a dual lamination pair $\Lambda^\pm_\phi \in \L^\pm(\phi)$, the following hold:
\begin{enumerate}
\item Any line $\ell \in \B$ that is not carried by $\A_\na(\Lambda^\pm_\phi)$ is weakly attracted either to $\Lambda^+_\phi$ by iteration of $\phi$ or to $\Lambda^-_\phi$ by iteration by $\phi^\inv$.
\item For any neighborhoods $V^+,V^- \subset \B$ of $\Lambda^+_\phi, \Lambda^-_\phi$, respectively, there exists an integer $m \ge 1$ such that for any line $\ell \in \B$ at least one of the following holds: \ $\gamma \in V^-$; \ $\phi^m(\gamma) \in V^+$; \ or $\gamma$ is carried by $\A_\na(\Lambda^\pm_\phi)$.
\end{enumerate}
\end{TheoremH}

\section*{Relatively irreducible subgroups (\PartFour).}

Our main result, Theorem~C, can be easily reduced to a special case, Theorem~I. For the rest of this section, when speaking about irreducibility relative to an extension of the form $\F \sqsubset \{[F_n]\}$, it is convenient to drop ``$\{[F_n]\}$'' from the terminology and speak of irreducibility relative to~$\F$.

\begin{TheoremI}
For each finitely generated subgroup $\h \subgroup \IA_n(\Z/3)$ and each $\h$-invariant free factor system $\F$, if $\F \sqsubset \{[F_n]\}$ is a multi-edge extension and if $\h$ is irreducible relative to $\F$ then then there exists $\phi \in \h$ which is fully irreducible relative to $\F$. 
\end{TheoremI}

The argument for reducing Theorem~C to Theorem~I is found in Section~\refRI{SectionReduction}: given an $\h$-invariant multi-edge extension $\F \sqsubset \F'$, one simply restricts the subgroup $\h$ to the stabilizer of the unique component of $\F'$ that is not a component of~$\F$, and then one applies Theorem~I to this restriction. 

Theorem~I is proved in Section~\refRI{SectionWholeShebang} using a ping pong argument which is developed in Section~\refRI{SectionFindingAttrLams}. The input for this argument is an outer automorphism $\psi \in \h$ and a lamination pair $\Lambda^\pm_\psi \in \L^\pm(\psi)$ having the property that $\Lambda^\pm_\psi$ is not supported by~$\F$ and hence $\F \sqsubset \A_\na(\Lambda^\pm_\psi)$; these are obtained by applying Theorem~E as discussed earlier. To simplify our exposition we shall assume that $\Lambda^\pm_\psi$ is nongeometric (but see Theorem~J for further discussion of this point). With this assumption, full irreducibility rel~$\F$ of $\psi$ is equivalent to the following properties of $\Lambda^\pm_\psi$: the nonattracting subgroup system $\A_\na\Lambda^\pm_\psi$ is equal to its minimal value~$\F$; and the joint free factor support $\F_\supp(\F,\Lambda^\pm_\psi)$ is equal to its maximal value $\{[F_n]\}$. The goal of the ping-pong argument is to take the given $\psi$ and $\Lambda^\pm_\psi$, for which these properties may fail, and produce a $\xi$ and $\Lambda^\pm_\xi$ for which these properties hold, and hence $\xi$ is fully irreducible rel~$\F$.

There are two ping-pong tournaments, each with several ping-pong games based on Proposition~\refRI{PropSmallerComplexity}. In each of these games one first forms some conjugate $\phi =\theta \psi \theta^\inv \in \h$ for some carefully chosen $\theta \in \h$, with corresponding lamination pair $\Lambda^\pm_\phi = \theta(\Lambda^\pm_\psi) \in \L^\pm(\phi)$. Then one takes a product of powers $\xi = \psi^l \phi^m$ for large exponents $l,m$, and one forms a certain lamination pair $\Lambda^\pm_\xi$ not supported by~$\F$.

In the first ping-pong tournament the goal is to drive down the value of $\A_\na\Lambda^\pm_\psi$ to its minimal value of $\F$. This is accomplished using one of the conclusions of Proposition~\refRI{PropSmallerComplexity}: if $l,m$ are sufficiently large then $\A_\na\Lambda^\pm_\xi$ is supported by each of $\A_\na\Lambda^\pm_\psi$ and $\A_\na\Lambda^\pm_\phi$, and hence is properly supported by each of them. By the chain condition on free factor systems, after repeated application one finds that $\A_\na\Lambda^\pm_\xi$ achieves the minimum value of~$\F$. 

In the second ping-pong tournament, restricting now to lamination pairs whose nonattracting subgroup system equals~$\F$ on the nose, the new goal is to drive up the value of $\F_\supp(\F,\Lambda^\pm_\psi)$ to its maximal value~$\{[F_n]\}$. This is accomplished using another of the conclusions of Proposition~\refRI{PropSmallerComplexity}: for sufficiently large $l,m$, the lamination $\Lambda^+_\xi$ is arbitrarily close to $\Lambda^+_\phi$ and the lamination $\Lambda^-_\xi$ is arbitrarily close to $\Lambda^-_\psi$. Further arguments using this conclusion show that $\F_\supp(\F,\Lambda^\pm_\xi)$ supports each of $\F_\supp(\F,\Lambda^\pm_\psi)$ and $\F_\supp(\F,\Lambda^\pm_\phi)$, and hence properly supports each of them. Again by the chain condition on free factor systems, after repeated application one finds that $\F_\supp(\F,\Lambda^\pm_\xi)$ achieves the maximum value of~$\{[F_n]\}$. 

The resulting $\xi$ is fully irreducible rel~$\F$.

\paragraph{Relatively geometric irreducible subgroups: Theorem J (Section~\refRI{SectionRelGeomIrr}).} In the context of Theorem~I, our discussions so far started with the assumed existence of $\psi \in \h$ and a lamination pair $\Lambda^\pm_\psi \in \L^\pm(\psi)$ that is not supported by~$\F$, such that $\Lambda^\pm_\psi$ is a nongeometric lamination pair. In fact, the reader may have noticed a bit of a swindle: not only did we assume that $\Lambda^\pm_\psi$ is nongeometric, but that all later lamination pairs encountered in the ping-pong tournament are nongeometric. In fact this swindle is justified by another conclusion of our main ping-pong results, Proposition~\refRI{PropSmallerComplexity}: assuming the input lamination pair $\Lambda^\pm_\phi$ is nongeometric, the output lamination pair $\Lambda^\pm_\xi$ produced by the ping-pong process will always be nongeometric as well.

Now let us consider the case that was set aside by this assumption: assume instead that the subgroup $\h$ is \emph{geometric relative to~$\F$} meaning that for all $\psi \in \h$ and all $\Lambda^\pm_\psi \in \L(\psi)$ not carried by~$\F$, the lamination pair $\Lambda^\pm_\psi$ is geometric. The ping-pong analysis can still be carried out in this case, and it gives interesting and much stronger conclusions. We state these conclusions here only in the absolute context where $\h$ is an irreducible subgroup that is absolutely geometric in that every attracting lamination of every element of $\h$ is geometric; the general, relative version of Theorem~J, stated in Section~\refRI{SectionRelGeomIrr}, requires the concepts of geometric models.

\begin{TheoremJ}[Absolute version]
Given a finitely generated subgroup $\h \subgroup \IA_n(\Z/3)$, if $\h$ is irreducible and geometric then there exists a compact surface $S$ with one boundary component and a homotopy equivalence between $S$ and the rose $R_n$ inducing an isomorphism $\Out(F_n) \approx \Out(\pi_1 S)$ under which $\h$ becomes a subgroup of the natural embedding $\MCG(S) \subgroup \Out(\pi_1 S)$, and $\h$ contains a pseudo-Anosov element of $\MCG(S)$. 
\end{TheoremJ}

The proof of Theorem~J starts out along the same lines as Theorem~I: one drives down $\A_\na(\Lambda_\psi)$ to its minimal value using a ping-pong argument (Section~\refRI{SectionProofUnivAttr}). The difference is that in the context where $\h$ is geometric rel~$\F$, once $\A_\na(\Lambda_\psi)$ has been driven down to its minimal value then the outer automorphism $\psi$ is already irreducible rel~$\F$; see Proposition~\refRI{PropUniversallyAttracting}. The step of ``driving up'' the free factor support, which is carried out in Section~\refRI{SectionLooking} for purposes of proving Theorems~C and~I, is not needed for Theorem~J. This is quite similar to the proof of the subgroup classification theorem for mapping class groups of surfaces carried out in \cite{Ivanov:subgroups}. 

The remaining steps of the proof of Theorem~J are carried out in Section~\refRI{SectionRelGeomIrr}. Since $\h$ is geometric irreducible rel~$\F$, the subgroup system $\A_\na(\Lambda_\psi)$ is \emph{never} a free factor system and in particular is never trivial. In the absolute case stated above, once $\A_\na(\Lambda_\psi)$ has been driven down to its minimum value it has a single infinite cyclic component corresponding to a one-holed surface $S$ as in the conclusion of Theorem~J. One uses minimality of $\A_\na(\Lambda_\psi)$ to show that this subgroup system is invariant under all of $\h$, from which the inclusion $\h \subgroup \MCG(S)$ follows. 

The general, relative case of Theorem~J proceeds similarly, the minimal value of $\A_\na(\Lambda_\psi)$ being the nonattracting subgroup system associated to an appropriate geometric model. See Section~\refRI{SectionRelGeomIrr} for the statement and proof of the relative case, and the remark following the statement for a reduction of the absolute case stated here to the general, relative case.

\bibliographystyle{amsalpha} 
\bibliography{mosher} 
 
 \printindex
 
 \end{document}